\documentclass[12pt]{amsart}
\usepackage{amsfonts}
\usepackage[pdftex]{graphicx}

\renewcommand{\epsilon}{\varepsilon}
\renewcommand{\phi}{\varphi}

\newtheorem{Lemma}{Lemma}
\newtheorem{Theorem}[Lemma]{Theorem}
\newtheorem{Proposition}[Lemma]{Proposition}

\newtheorem{Definition}[Lemma]{Definition}

\begin{document}

\address{Department of Mathematics, University of Pittsburgh, 301 Thackeray Hall, Pittsburgh, PA 15260, USA.}
\author{Alexander Borisov}
\thanks{The research of the author was supported in part by the NSA grants H98230-08-1-0129 and H98230-06-1-0034}
\title[Jacobian Conjecture]{A geometric approach to the two-dimensional Jacobian Conjecture}
\email{borisov@pitt.edu}
\maketitle

\section{Introduction}

Suppose $f(x,y), g(x,y)$ are two polynomials with complex coefficients. The classical Jacobian Conjecture (due to Keller) asserts the following.

{\bf Conjecture.} (Jacobian Conjecture in dimension two) If the Jacobian of the pair $(f,g)$ is a non-zero constant, then the map $(x,y)\mapsto (f(x,y),g(x,y))$ is invertible. Note that the opposite is clearly true, because the Jacobian of any polynomial map is a polynomial, and, when the map is invertible, it must have no zeroes, so it is a constant.

The Jacobian Conjecture and its generalizations received considerable attention in the past, see \cite{Essen}. It is notorious for its subtlety, having produced a substantial number of wrong "proofs", by respectable mathematicians. 

From the point of view of a birational geometer, the most natural approach to the two-dimensional Jacobian Conjecture is the following. Suppose a counterexample exists. It gives a rational map from $P^2$ to $P^2.$ After a sequence of blowups of points, we can get a surface $X$ with two maps: $\pi : X \to P^2$ (projection onto the origin $P^2$) and $\phi : X \to P^2$ (the lift of an original rational map).

Note that $X$ contains a Zariski open subset isomorphic to $A^2$ and its complement, $\pi ^* ((\infty))$, is a tree of smooth rational curves. We will call these curves exceptional, or curves at infinity. The structure of this tree is easy to understand inductively, as it is built from a single curve $(\infty)$ on $P^2$ by a sequence of two operations: blowing up a point on one of the curves of blowing up a point of intersection of two curves. However, a non-inductive description is probably impossible, which is the first difficulty in this approach. Another difficulty comes from the fact that the exceptional curves on $X$ may behave very differently with respect to the map $\phi .$ More precisely, there are four types of curves $E$.

type 1) $\phi (E) = (\infty )$
 
type 2) $\phi (E)$ is a point on $(\infty )$

type 3) $\phi (E)$ is a curve, different from $(\infty)$ 

type 4) $\phi (E)$ is a point not on $(\infty)$

From a first glance, the situation appears almost hopeless. The goal of this paper is to show that it really is not that bad. In particular, for a given graph of curves, one can essentially always tell which curves are of which type, and there is a fairly restrictive family of graphs that can potentially appear in a counterexample to the JC.

{\bf Acknowledgments} The author is indebted to V.A. Iskovskikh for introducing him to the beauty of birational geometry, and to David Wright and Ed Formanek for stimulating discussions related to the Jacobian Conjecture.  

\section{Preliminary observations}
We change the notation slightly from the Introduction.

Suppose $X= P^2$, $Y=P^2$ and $\phi _X^Y: X ---> Y $ is a rational map. Suppose further that on an open subset $A^2\subset P^2 =X$ the map $\phi$ is defined, unramified, and $\phi (A^2) \subseteq A^2 \subset P^2 =Y.$ By a sequence of blowups at smooth points, we get a surface $X_1$ with a birational map $\pi : X_1 \to X $ and a generically finite map $\phi = \phi_{X_1}^Y : X_1 \to Y$ such that $\phi _{X_1}^Y = \phi _X^Y \circ \pi .$
 
The blowups that lead to $X_1$ can be done outside of $A^2\subset X.$ So $X_1 = A^2 \cup (\cup E_i)$, where $E_i$ are rational curves. The following observation is straightforward. 

\begin{Proposition}

1) The curves $E_i$ form  a tree.

2) One of $E_i$ is $\pi ^{-1} (\infty),$ all others are mapped to points by $\pi $.

3) The classes of $E_i$ form a basis in the Picard group of $X_1$.

\end{Proposition}

The structure of $X_1$ is largely determined by the graph of intersections of $E_i$. The vertices of this graph correspond to $E_i$-s and are usually labeled by $(-E_i^2).$ The edges correspond to the points of intersections of two different $E_i$-s. The graph is a tree.

This graph is not so easy to deal with because blowing up a point changes the self-intersections of the curves passing through it. Inspired by the Minimal Model Program, we consider a different labeling of this graph. We consider the augmented canonical class of $X_1,$ $\bar{K}_{X_1} = K_{X_1} +\sum _i E_i$. It can be uniquely written as a linear combination of $E_i,$ $\bar{K}_{X_1}= \sum _i a_iE_i$. We label the vertices of the intersection graph by these numbers $a_i.$

With this labeling we now describe what happens when a point is blown up, in any of the  intermediate steps in getting from $X$ to $X_1.$ 

\begin{Proposition}
When a point is blown up, going from $X_1'$ to $X_1'',$ one of the following two operations is performed to the graph of exceptional curves:

1) A new vertex is added to the graph, connected to one of the vertices. It is labeled $a_i+1,$ where $a_i$ is the label of the vertex it is connected to.

2) A new vertex is introduced on the edge connecting two vertices, "breaking" the edge into two edges. The new vertex gets labeled with $a_i+a_j,$ where $a_i$ and $a_j$ are the labels of the two vertices it is connected to.
\end{Proposition}

{\bf Proof.} The first case corresponds to blowing up a point on one of the curves. The second case corresponds to blowing up an intersection of two curves. The augmented canonical class calculations are straightforward and are left to the reader.

Notice that once a vertex is created, its label never changes, which is in sharp contrast with the traditional labeling.

The following observation is true for any $X_1,$ unrelated to the map $\phi .$  It is easily proven by induction on the number of exceptional curves, using the above proposition.

\begin{Proposition} For any two adjacent vertices $E_i,$ $E_j$ of the graph of $X_1,$ $gcd(a_i,a_j)=1.$ In particular, no two adjacent vertices have even labels.
\end{Proposition}

The following example serves two purposes. It shows how the graph of $X_1$ is constructed from the graph of $X=P^2,$ and we will use it to compare our labeling with the traditional self-intersection labeling.

{\bf Example.}  We start with $X=P^2,$ its graph is the following.

$$\circ$$
\vskip -0.6cm
$$-2$$

Blowing up a point, we get
$$\circ \!\! - \!\!\! - \!\!\! - \!\!\! - \!\! \circ$$
\vskip -0.6cm
$$-1\ \  -2$$

Blowing up another point, on the pullback of $(\infty),$ we get 

$$\circ \!\! - \!\!\! - \!\!\! - \!\!\! - \!\! \circ \!\! - \!\!\! - \!\!\! - \!\!\! - \!\! \circ $$
\vskip -0.6cm
$$-1\ \  -2 \ \ -1$$

Blowing up a point on a newly blown up curve, we get 

\newpage

$$\circ \!\! - \!\!\! - \!\!\! - \!\!\! - \!\! \circ \!\! - \!\!\! - \!\!\! - \!\!\! - \!\! \circ \!\! - \!\!\! - \!\!\! - \!\!\! - \!\! \circ  $$
\vskip -0.6cm
$$-1\ \  -2 \ \ -1 \ \ \ \  0$$

Then we blow up the intersection of the last two curves and get

$$\circ \!\! - \!\!\! - \!\!\! - \!\!\! - \!\! \circ \!\! - \!\!\! - \!\!\! - \!\!\! - \!\! \circ \!\! - \!\!\! - \!\!\! - \!\!\! - \!\! \circ  \!\! - \!\!\! - \!\!\! - \!\!\! - \!\! \circ  $$
\vskip -0.6cm
$$-1 \ \  -2 \ \ -1 \ \ \ -1 \ \ \ \ \  0\ $$

Blowing up another intersection point, we get 

$$\circ \!\! - \!\!\! - \!\!\! - \!\!\! - \!\! \circ \!\! - \!\!\! - \!\!\! - \!\!\! - \!\! \circ \!\! - \!\!\! - \!\!\! - \!\!\! - \!\! \circ  \!\! - \!\!\! - \!\!\! - \!\!\! - \!\! \circ  \!\! - \!\!\! - \!\!\! - \!\!\! - \!\! \circ  $$
\vskip -0.6cm
$$-1 \ \  -2 \ \ -1 \ \ \ -2 \ \ \  -1 \ \ \ \ \  0\ $$

Blowing up another point, we get the following graph

$$\circ \!\! - \!\!\! - \!\!\! - \!\!\! - \!\! \circ \!\! - \!\!\! - \!\!\! - \!\!\! - \!\! \circ \!\! - \!\!\! - \!\!\! - \!\!\! - \!\! \circ  \!\! - \!\!\! - \!\!\! - \!\!\! - \!\! \circ  \!\! - \!\!\! - \!\!\! - \!\!\! - \!\! \circ \!\! - \!\!\! - \!\!\! - \!\!\! - \!\! \circ  $$
\vskip -0.6cm
$$\ 0\ \ \ -1 \ \  -2 \ \ \ -1 \ \ \ -2 \ \ \  -1 \ \ \ \ \  0\ $$

Finally, blowing up four more points (in any order) we get the following:

\hskip 2.37cm $1\circ$ \hskip 6.15cm $\circ 1$

\vskip -0.2cm

\hskip 2.6cm $|$ \hskip 6.2cm $|$

\vskip -0.2cm

\hskip 2.6cm $|$ \hskip 6.2cm $|$

\vskip -0.7cm

$$\circ \!\! - \!\!\! - \!\!\! - \!\!\! - \!\! \circ \!\! - \!\!\! - \!\!\! - \!\!\! - \!\! \circ \!\! - \!\!\! - \!\!\! - \!\!\! - \!\! \circ  \!\! - \!\!\! - \!\!\! - \!\!\! - \!\! \circ  \!\! - \!\!\! - \!\!\! - \!\!\! - \!\! \circ \!\! - \!\!\! - \!\!\! - \!\!\! - \!\! \circ  $$
\vskip -0.7cm
$$0\ \ -1 \ \  -2 \ \ -1 \ \  -2 \ \  -1 \ \ \ \  0$$

\vskip -0.65cm

\hskip 2.6cm $|$ \hskip 6.2cm $|$

\vskip -0.2cm

\hskip 2.6cm $|$ \hskip 6.2cm $|$

\vskip -0.2cm

\hskip 2.37cm $1\circ$ \hskip 6.15cm $\circ 1$

For most of the exceptional curves, one can easily recover their self-intersection from the graph, using the adjunction formula:

$(K_{X_1}+E_i)E_i=-2,$ so $\bar{K}_{X_1}\cdot E_i=-2+\# (E_j \ {\textrm{adjacent to }}E_i)$

Thus, if $\bar{K}_{X_1} = \sum a_i E_i,$ we have 
$$a_iE_i^2 + \sum \limits_{E_j \ adj. \ E_i} a_j = -2 + \# (E_j \ {\textrm{adjacent to }}E_i)$$

So if $a_i\neq 0,$ $E_i^2$ can be easily calculated.

However, when $a_i=0,$ it is not that easy. One can see in the above example, the left curve with $a_i=0$ has self-intersection $(-3),$ while the right one has self-intersection $(-4),$ despite the symmetry of the graph. One can remedy this situation by keeping track of the strict pullback of infinity. We do not need it in this paper, and the details are left to an ambitious reader.

Note that the subgraph of vertices with negative labels is connected. It is separated from the "positive" vertices by the "zero" vertices. Moreover, the "zero" vertices are only connected to vertices with labels $(-1)$ or $1$.

Now we are going to make use of the map $\phi .$ The main idea is to use the adjunction formula for $\phi$ to get a formula for $\bar{K}_{X_1}.$

Recall from the Introduction the four types of curves $E_i.$ For every curve of type $1$ or $3$ denote by $d_i$ the degree of its image on $Y=P^2,$ by $f_i$ the degree of the map onto its image and by $r_i$ the ramification index. Denote by $L$ the class of the line on $Y=P^2.$

\begin{Proposition} There exist integers $b_i$ for the curves $E_i$ of types $2$ and $4$ such that
$$\bar{K}_{X_1} = \phi ^* (-2L) + \sum \limits_{type(E_i)=3} r_iE_i + \sum \limits_{type(E_i) =2 or 4} b_iE_i$$
\end{Proposition}

{\bf Proof.} Consider the differential form $\omega $ on $Y=P^2$ that has the pole of order $3$ at $(\infty)$ and no other poles or zeroes. Because $\phi $ is unramified on the $A^2 \subset X,$ there is a differential form on $X_1,$ such that its divisor of zeroes and poles is $\phi ^*(-3L) + \sum _i c_iE_i,$ where $c_i$ can be calculated locally at a general point of each $E_i.$ 

Notice that for the curves $E_i$ of types $1$ and $3$, $c_i=r_i-1,$ and 
$$\phi ^* (L) = \sum \limits_{type (E_i)=1} r_iE_i + \sum \limits_{type (E_i)=2} e_iE_i$$
for some $e_i$. Thus,
$$\bar{K}_{X_1} = K_{X_1} + \sum E_i = \phi ^* (-3L) + \sum \limits_{type(E_i)=1 or 3} r_iE_i + \sum \limits_{type(E_i)=2 or 4} (c_i+1)E_i=$$
$$=\phi ^*(-2L) + \sum \limits_{type (E_i)=3} +\sum \limits_{type (E_i) =2 or 4} b_i E_i$$
Q.E.D.

Note that because $E_i$ are independent in the Picard group of $X_1,$ the above representation of $\bar{K}_{X_1}$ is unique and must match with the labeling of the graph of $E_i$. As a corollary, we have the following observation.

\begin{Proposition}
1) Any curve of type $1$ has a negative even label.

2) Any curve of type $3$ has a positive label.
\end{Proposition}

{\bf Proof.} Note that $\phi ^*(-2L)$ only involves curves of type $1$ and $2$. Q.E.D.
 
Additionally, the union of curves of type $1$ and $2$ must be connected, as a specialization (set-theoretically) of a pullback of a generic $L$ on $Y=P^2.$ This means that the corresponding subgraph is connected. 

Every curve of type $3$ must intersect with one of the curves of type $1$ or $2$, while the curves of type $4$ do not intersect with curves of type $1$ or $2$. (This follows from the projection formula of the intersection theory: if $E$ is  a curve on $X_1,$ $E\cdot \phi^* (L) = (\phi_* E) \cdot L.$)

On the other hand, a type $3$ curve can not intersect a type $1$ curve, because negative and positive labels are never adjacent. Because the graph of the exceptional curves on $X_1$ is a tree, no two curves of type $3$ intersect with each other. Putting this all together, we must have the following. The tree of curves on $X_1$ has a connected subtree containing all curves of type $1$ and $2$. Some of the vertices of this subtree may have one or more curves of type $3$ connected to them. Then some of these type $3$ curves may have trees of type $4$ curves connected to them. Additionally, no two curves of type $1$ are adjacent, and the subtree of curves of type $1$ and $2$ contains the connected subtree of curves with negative labels.

\begin{Proposition} $\pi ^{-1} (\infty )$ is of type $1$ or $2$.
\end{Proposition}

{\bf Proof.} One can prove it using the above description of the graph of exceptional curves, but there is also the following direct geometric argument. The pullbacks of lines on $X=P^2$ form a family of rational curves $C$ on $X_1$ that intersect $\pi ^{-1} (\infty)$ at one point while being inside $A^2$ elsewhere. Consider $\phi (C)$ for a generic $C.$ If $\pi ^{-1} (\infty) $ is of type $3$ or $4$ then $\phi (C) \subseteq A^2 \subset Y.$ The curve $C$ is proper and $A^2$ is affine, so $\phi (C)$ is a point, which is impossible. Q.E.D.

Until now, the variety $X_1$ was an arbitrary resolution at infinity of the original rational map. But we can put an additional restriction on it, to avoid unnecessary blowups.

\begin{Definition}
A curve $E_i$ on $X_1$ is called {\bf final} if there is a sequence of blowups from $X$ to $X_1$ such that $E_i$ is blown up last.
\end{Definition}

Note that there may be more than one final curve, and $\pi ^{-1} (\infty)$ is never final. In what follows, $E_i$ is one of the exceptional curves on $X_1$.

\begin{Proposition}
Suppose that when $X_1$ was created, $E_i$ was created after all of its neighbors in the graph (i.e. all adjacent vertices). Then $E_i$ is a final curve.
\end{Proposition}

{\bf Proof.} Instead of creating $E_i$ at its due time we can change the order of blowups and create it at the last step of the process, without changing anything else. Q.E.D.
  
\begin{Proposition}
Suppose $a_i=a(E_i)\geq 2$ and it is the maximum (not necessarily a strict maximum) label among all its neighbors. Then $E_i$ is final.
\end{Proposition}

{\bf Proof.} We will prove that $E_i$ was created after all its neighbors. First of all, no neighbor of $E_i$ can be a blowup of a point on $E_i,$ because its label would have been $a_i+1.$ If it were a blowup of a point of intersection of $E_i$ and some $E_j,$ then before the blowup we had $E_i$ and $E_j$ were adjacent before the blowup. Negative curves are never adjacent to the positive curves and zero curves are only adjacent to curves with labels $1$ of $-1.$ Thus, $a_j \geq 1.$ So the label of the new curve is $a_i+a_j\geq a_i+1>a_i.$ Q.E.D.

Note that no two curves with the same label $a_i\geq 2$ can be adjacent, by Proposition 3. So every local maximum $a_i\geq 2$ is a strict maximum.

\begin{Proposition}
If $a_i=1,$ then $E_i$ is final if and only if it either has only one neighbor, with label $0$, or exactly two neighbors, with labels $1$ and $0.$ 
\end{Proposition}

{\bf Proof.} A curve with label $1$ can be created either by a blowup of a point on a curve with label $0$ or by a blowup of an intersection of a curve with label $0$ and a curve with label $1$. Once created, it will be final if and only if no other curve is blown up as its neighbor. The rest is easy and is left to the reader. 

The above two propositions allow us to easily spot the final curves in the positive part of the graph of curves. Our interest in the final curves stems from the following. If one of the final curves on $X_1$ is of type $2$ or $4$, then it can be contracted, using the $\phi -$relative MMP to get another $X_1,$ with two maps to $X$ and $Y$ and a smaller Picard number.

\begin{Definition}
We call $X_1$ {\bf minimal} if all of its final curves are of type $1$ or $3$.
\end{Definition} 

\begin{Proposition}
If a counterexample to JC exists, it can be obtained using a minimal $X_1.$
\end{Proposition}

{\bf Proof.} Just look at the $X_1$ with smallest possible Picard number. If it is not minimal, it can be created in such a way so that some curve of type $2$ or $4$ is blown-up last. Using MMP relative to $\phi,$ it can be blown down, maintaining the morphisms, and creating a counterexample to JC with smaller Picard number. Q.E.D.

From now on, $X_1$ will always be minimal.

\begin{Proposition} Suppose $E$ is a curve of type $3$ on $X_1.$ Suppose $E_0$ is the curve of type $2$ it is adjacent to. Then the tree on the other side of $E$ is a line $E-E_1-...-E_k,$ where $E_1,...E_k$ are of type $4$. Additionally, the order of creation is the following: $E_0,E_k,E_{k-1},...,E_1,E.$
\end{Proposition}

{\bf Proof.} The label of $E$ is positive and all curves $E,$ $E_1,...E_k$ must be created after the curve $E_0$. The last one created must be of type $3,$ for any possible order of creation, from which the result follows. Q.E.D.

\section{Other varieties and further analysis}

We start with the theorem that shows that type $3$ curves must exist in a counterexample to the JC. This result is well known and can be easily proven by a topological argument, so the main purpose of our proof is to show an easy application of our method before proceeding to the more intricate questions.

\begin{Theorem} Suppose $X_1$ and $\phi$ provide a counterexample to the JC. Then $X_1$ contains a curve of type $3$, where $\phi$ is ramified. 
\end{Theorem}

{\bf Proof.} 
Consider a generic line $L$ on the target variety $Y=P^2.$ The curve $C=\phi^{-1}(L)$ is smooth and irreducible ("Bertini's theorem"). Moreover, we can assume that for all but finitely many lines $L'$ that only intersect $L$ "at infinity",  $C'=\phi^{-1}(L')$ is smooth and irreducible. We can also assume that $L$ does not pass through the images of the exceptional curves of types $2$ and $4,$ so $C$ does not intersect these curves on $X.$ Suppose that the genus of $C$ is $g$, the map $H=\phi_{|_{C}}C\to L $ has degree $n$ and the number of points of $C$ "at infinity" is $k$. (There is a special point $\infty$ on $L$, the only one not lying in $A^2.$ The number $k$ is the number of points of $C$ mapped to it, in a  set-theoretic sense.) Because the map $\phi$ is only ramified at the exceptional curves of $X,$ the map $H$ could only by ramified at these $k$ points at infinity. By Hurwitz formula, we have 
$$2g-2=-2n+r,$$
where $r$ is the total ramification at infinity. We have $g\geq 0,$ $n\geq 1$ and $r\leq n-k.$ So

$$-2\leq 2g-2 \leq -2n+n-k = -n-k \leq -2 $$
Thus all the inequalities above are equalities, $g=0,$ $n=1,$ and $k=1$. 

The curve $\pi (C)$ is now a rational curve on $P^2,$ and its part on $A^2$ is smooth, while its intersection with the infinity consists of one point. By the celebrated Abhyankar-Moh-Suzuki Theorem (see \cite{Essen}) there is a polynomial automorphism of $A^2$ that maps it to the standard line $y=0.$ Precomposing the map $\phi$ with this automorphism, we now assume that the line $y=0$ is mapped by $\phi$ to a line on $A^2.$ By a linear transformation in the target variety we assume that it is mapped to the line $y=0$ there.

The map $\phi$ restricted to this line is $1-to-1$ so it is of the form $\phi (x,0)=(ax+b,0)$. By another linear transformation we can assume that  it is $\phi (x,0)=(x,0).$ Now consider the preimages by $\phi$ of the lines $y=c$ on the target variety, for generic $c$. By the same argument as before these curves are rational and smooth. They do not intersect $y=0$ except "at infinity". Because these curves are also isomorphic to $A^1,$ they are given by two polynomials $(x(t),y(t))$. Because $y(t)$ is never zero, it must be a constant. So these curves are of the form $y=d$ for some constant $d$ that depends on $c$. Thus the second coordinate function of the map $\phi$ only depends on $y$. This is a polynomial and it must be a linear one. By a linear change of variables, we can assume that it is just $y$. Now the map $\phi$ has the following form:
$$\phi(x,y)=(a(y)x+b(y),y)$$
Clearly, the polynomial $a(x)$ must be a constant. The map $\phi$ is then easily invertible, thus $X_1$ is not a counterexample to the JC. Q.E.D.

Now we want to make further use of the morphism $\phi : X_1 \to Y.$ First, we decompose it into a composition of two morphisms, birational and finite (Stein decomposition):
$$X_1 \longrightarrow Y_1 \longrightarrow Y$$ 

Here the first morphism is birational and denoted $\phi_1$, and the second one is finite and denoted $\phi_2.$ 

Here $Y_1$ is a normal surface, and one can talk about its canonical class, defined modulo numerical equivalence. By adjunction, we have:
$$K_{Y_1} = \phi_2^* K_Y + \sum (r_i-1) E_i,$$
where $r_i$ is the ramification index, and $E_i$ are dimension $1$ images of the curves $E_i$ of types $1$ and $3$ on $X_1.$

Define $\bar{K}_[Y_1]=K_{Y_1}+\sum E_i$.
Then $\bar{K}_{Y_1}$ looks very simple:
$$\bar{K}_{Y_1}=\phi_2^*(-3L)+\sum r_iE_i = \phi_2^*(-2L) +\sum \limits_{type(E_i)=3} r_iE_i.$$

Additionally, we decompose $\phi _1$ using the Log Minimal Model Program for $(X_1, \sum E_i)$ relative to $\phi_1:$
$$X_1\longrightarrow X_2 \longrightarrow Y_1$$
Here the first map is $\phi_3$, and the second one is $\phi_4.$ Note hat $(X_1, \phi_3( \sum E_i))$ has  log-terminal singularities and $\bar{K}_{X_2} = K_{X_2} + \sum \phi_3 (E_i)$ is $\phi_4 -$nef.  

Because $\bar{K}_{X_2}$ is $\phi_4-$nef, 
$$\bar{K}_{X_2} = \phi_4^* (\bar{K}_{Y_1}) + \sum \limits_{\phi_{4}(E_i)=pt} c_i E_i, \ \ {\textrm{where}} \ \ c_i\leq 0$$
So 
$$\bar{K}_{X_2} = \phi_4^*\phi_2^* (-2L) + \phi_4^*(\sum \limits_{type(E_i)=3} r_iE_i) + \sum \limits_{\phi_{4}(E_i)=pt} c_i E_i $$

The following calculation is very important for the paper.
Consider $\bar{K}_{X_2} \cdot (\sum \limits_{type(E_i)=3} r_iE_i) .$ It equals
$$-2\sum \limits_{type(E_i)=3} (r_id_if_i) + \sum \limits_{\phi_4(E_j)=pt}  c_j \sum \limits_{type(E_i)=3} r_iE_i + \phi_4^* (\sum \limits_{type(E_i)=3} r_iE_i) \sum \limits_{type(E_i)=3} r_iE_i$$

On the other hand, 
$$\bar{K}_{X_2} \cdot E_i = (K_{X_2}+E_i)E_i +\sum \limits_{E_j\neq E_i} E_jE_i \geq -2 +1 =-1 $$

So, $\bar{K}_{X_2} \cdot \sum r_iE_i \geq - \sum r_i.$ Thus,
$$\sum r_i(-2d_if_i) + \phi_4^* (\sum \limits_{E_i\subset Y_1, type(E_i)=3} r_iE_i) \cdot \sum \limits_{E_i\subset X_2, type(E_i) =3} r_iE_i \geq -\sum r_i$$

Note that
$$\phi_4^*(\sum r_iE_i)\cdot \sum r_iE_i =(\phi_4^* \sum \limits_{E_i\subset Y_1, type(E_i) = 3} r_iE_i)^2$$

So we have

$$(\phi_4^* \sum \limits_{E_i\subset Y_1, type(E_i) = 3} r_iE_i)^2 \geq \sum r_i(2d_if_i -1 ) \geq 0.$$

Pulling back to $X_1,$ we get 
$$(\sum\limits_{E_i\subset Y_1, type (E_i) =3} r_i \phi_1^* (E_i))^2 \geq 0$$

\begin{Proposition} The curve $\pi ^{-1} (\infty)$ is of type $2.$
\end{Proposition}

{\bf Proof.} By Proposition 6, it is of type $1$ or $2$. If it is of type $1$, then it is not included in $\phi_1^* \sum r_iE_i$ above. So $\phi_1^*(\sum r_iE_i)$ above consists of curves contractible by $\pi .$ So, its self-intersection is negative, contradicting the calculation above. Q.E.D.

Now we define a graph of curves on $X_2,$ similar to that on $X_1.$ Notice that $(X_2, \sum E_i)$ has log terminal singularities. As a result, the curves $E_i$ have simple normal crossings at smooth points of $X_1.$ The only difference from $X_1$ is that the curves may pass through singularities of $X_2.$ Denote this graph by $\Gamma _{X_2}.$ It is obtained from $\Gamma_{X_1}$ by contracting all type $4$ curves and some type $2$ curves. 
Suppose $E_i \subset X_2$ is of type $3$.

Note that every curve of type $3$ on $X_1$ intersects the union of curves of type $2$ at exactly one point, and does not intersect curves of type $1.$ When the curves of type $2$ are contracted, on $Y_1,$ every curve of type $3$ intersects the union of curves of type $1$ at exactly one point.

\begin{Proposition}
For every curve $E_i$ of type $3$ on $Y_1$ the point above is $\phi_1(\pi^{-1} (\infty))$.
\end{Proposition} 

{\bf Proof.} Suppose there is a point $y\in Y_1$ on the union of type $1$ curves, which is not $\phi_1 (\pi^* (\infty))$ and which has some type $3$ curves passing through it. Consider
$$\bar{K}_{Y_1} \cdot \sum \limits_{y\in E_i, type(E_i)=3} r_iE_i$$

On the one hand,
$$\bar{K}_{Y_1} \cdot E_i \geq (K_{Y_1}+E_i)E_i \geq -2,$$

so 

$$\bar{K}_{Y_1} \cdot \sum \limits_{y\in E_i, type(E_i)=3} r_iE_i \geq -2 \sum  \limits_{y\in E_i, type(E_i)=3} r_i$$

On the other hand,

$$\bar{K}_{Y_1} \cdot E_i \geq (K_{Y_1}+E_i)E_i  =-2 \sum \limits_{y\in E_i, type (E_i)=3} r_iE_i = $$

$$=-2 \sum \limits_{y\in E_i, type (E_i)=3} (f_id_i) + \left( \sum \limits_{y\in E_i, type (E_i)=3} r_iE_i \right)^2$$

(Note that curves of type $3$ can only intersect at the union of curves of type $1.$)

So,

$$ \left( \sum \limits_{y\in E_i, type(E_i)=3} r_iE_i \right)^2 =2\sum \limits_{y\in E_i, type(E_i)=3}  r_i (f_id_i-1) \geq 0.$$

But, pulling back to $X_1,$ 
$$ \left( \sum \limits_{y\in E_i, type(E_i)=3} \phi_4^*(r_iE_i) \right)^2 <0,$$

because $ \sum \limits_{y\in E_i, type(E_i)=3} \phi_4^*(r_iE_i)$ does not include $\pi ^{-1}(\infty),$ contradiction. Q.E.D.

\begin{Proposition}
On $Y_1,$ all exceptional curves contain $\phi_1 (\pi ^{-1} (\infty))$ and there are no other points of intersection.
\end{Proposition}

{\bf Proof.} By the proposition above, every curve of type $3$ contains $\phi_1 (\pi ^{-1}(\infty))$ and this is the only point of intersection with other curves. Now consider a curve $E_i$ of type $1$. Suppose it does not contain $\phi_1 (\pi^{-1}(\infty))$. Then it does not intersect any of the curves of type $3$ on $Y_1.$ 

On $Y_1$ we have:
$$\bar{K}_{Y_1} \cdot E_i \geq  (K_{Y_1}+E_i) E_i \geq -2 $$

On the other hand,
$$\bar{K}_{Y_1} \cdot E_i = (-2\phi_2^*(L)+\sum \limits_{type (E_j) =3} r_j E_j ) \cdot E_i = -2\phi_2^*(L) \cdot E_i \leq -2 $$

The inequalities above became equalities if and only if $E_i$ intersects no other curves and is smooth. This would make it the only curve of type $1$, which would have to intersect with some curves of type $3,$ contradiction. Q.E.D.

Thus, we know that every curve of type $1$ on $Y_1$ contains $\phi_1 (\pi ^{-1} (\infty)).$ We now look at the graph of curves on $X_1.$ The curves of type $2$ that are mapped to $\phi_1 (\pi ^{-1} (\infty))$ form a connected subgraph, containing $\pi ^{-1} (\infty).$ Every curve of type $1$ or $3$ is attached to this subgraph. On "the other side" of the curves of type $3$ there may be curves of type $4$, and on "the other side" of the curves of type $1$ there may be curves of type $2$. Note that all of these "other side" curves must be created before the corresponding type $3$ or type $1$ curves. When mapped to $Y_1,$ the curves of type $1$ and $3$ intersect at $\phi_1 (\pi ^{-1} (\infty))$ and nowhere else.

One can restrict the structure of possible counterexamples even further.
\begin{Theorem} In any counterexample to the Jacobian Conjecture there are no curves on "the other side" of the curves of type $1$.
\end{Theorem}

{\bf Proof.} Consider a curve of type $1$, $E,$ on $X_1$. Suppose the ramification index at $E$ is $r$. Then the coefficient of $\phi ^* (L)$ in $E$ is $r$, and the coefficient of $\bar{K_{X_1}}$ is $(-2r).$ Consider the divisor class $D=\bar{K_{X_1}}+2\phi* (L)= ...+0\cdot E +x_1E_1+...+x_kE_k,$ where $E_1,...E_k$ are the curves on $X_1$ "on the other side" of $E$.
We know that $D$ intersects by zero with $E_1,...,E_{k-1}$. It intersects by $-1$ with $E_k$. We formally add another vertex to the graph, $"E_{k+1}"$ and set the coefficient of $D$ at it to be $1$. (Note that we are not blowing up any points and $E_{k+1}$ does not have any geometric meaning. We now have a chain $E, E_1, ..., E_{k}, E_{k+1}$ and a divisor $D'= 0\cdot  E +x_1E_1+...+x_kE_k +1\cdot E_{k+1},$ such that $D'$ intersects by zero with all $E_1, E_2,...,E_k.$ Because the self-intesections of all $E_i,$ $1\leq i \leq k,$ is less than or equal to $-2,$ the coefficients $x_i$ must form a concave up chain between $0$ and $1$, contradicting their integrality. (Here is a more formal argument. Suppose at least one of the $x_i,$ $1\leq i \leq k,$ is not positive. Then consider the minimum of $x_i$, obtained at $x_j$, such that $x_{j+1}> x_j$, where formally $x_0=0, x_{k+1}=1$. Then $D'\cdot E_j \geq x_{j-1}+x_{j+1} +2x_j >0 ,$ contradiction.  Suppose the maximum of $x_i,$ $1\leq i \leq k,$ is greater than or equal $1$ and is obtained at $x_j,$ where $x_{j-1}<x_j.$  Then $D'\cdot E_j leq x_{j-1}+x_{j+1} - 2 x_j <0,$ contradiction. Thus all $x_i$ are strictly between $0$ and $1$, which is impossible because they are integers.) Q.E.D.

One can restrict the graph of curves on $X_1$ even further, by using yet another labeling on the exceptional curves. We denote by $Q_{X_1}$ the anti-self-intersection matrix, which is minus the matrix of the graph. The following lemma is easy to check.

\begin{Lemma} Suppose $E$ is one of the exceptional curves on $X_1.$ Denote by $d_E$ the determinant of the submatrix of the anti-selfintersection matrix $Q_{X_1}$, obtained by removing $E$. Then any additional blowups of $X_1$ do not change $d_E$. 
\end{Lemma}

{\bf Proof.}  This submatrix is minus the matrix of the graph, which is obtained from the graph of $X_1$ by simply removing $E$ and all edges at $E$. There are two kinds of blowups to consider: a blowup of a new point on $E$ and a blowup of the intersection of $E$ and another exceptional curve, $E_1$. In the former case, the new minor corresponds to the graph which is obtained from the graph of the original minor by adding extra vertex, with self-intersection $(-1),$ not connected to anything else. The corresponding determinant is clearly $1$ times the original minor. In the latter case, the new submatrix is obtained by adding a new vetex, $E_2$, with the diagonal entry $1,$ and two entries $(-1)$ that correspond to the edges with $E_1$. Also, the diagonal entry of $E_1$ is increased by $1$. Using cofactor expansion at $E_2$ (and then at $E_1$), one can easily check that the determinant does not change. Q.E.D.

We will call the new labeling described in the above lemma as the determinant labeling. Its significance follows from the fact that the determinant label of any exceptional curve of type $1$ must be negative. Otherwise, we get a contradiction, like in Proposition 15.

The determinant labeling is not very easy to control. When a new point is blown up, its label is one plus the label of the  "parent" curve.  However the label of the blowup of the intersection of two curves depends on more than just the labels of these curves. One can prove however that for an exceptonal curve of type $1$ to have a negative determinant label, at least one of its "ancestors" (curves that must be blown up before it) must have $\bar{K}$-label zero. We will not directly use this fact, and its proof is omitted from this preliminary version of the paper. This was however used as a guide to help find the Example in the next section.

\section{The search for the counterexamples}

We have now established that the structure of the graph of exceptional curves on a possible counterexample to the Jacobian Conjecture is very far from arbitrary. Besides giving hope for the proof of the conjecture, this also greatly helps to look for a counterexample to it, if it exists. Every such counterexample must be given by a linear system on some surface $X_1$ with very special properties. Suppose this linear system is a subsystem of the complete system $|L|$, $L=\sum {a_i}E_i$.

1) For every final curve with negative $\bar{K}$-label, this label must be even, and the coefficient $a_i$ is minus half of it.

2) The intersection of $L$ with each of the above curves is $1$. 

3) The numbers $a_i$ are zero for all final curves with positive $\bar{K}$-label, and all curves "on the other side" of them. 

4) $L$ intersects by zero with all non-final curves connected to the pullback of infinity by non-final curves.

Additionally, every final curve with negative $\bar{K}$-label (a hypothetical curve of type $1$) must contain a $0$-curve in the set of its ancestors. Playing with different possibilities for the blown-up curves, one can produce a number of pairs $(X_1, \{a_i\})$ that satisfy the above condition. For each of these, one can easily calculate a lower bound on $h^0(L)$ using the Riemann-Roch formula. (Note that $h^2(L)=h^0(K-L)$ is zero, because any divisor, which is rationally equivalent to $K-L$, must intersect negatively with the pull-back of a line from the original $P^2,$ so can not be effective). In many cases this estimate, $\frac{L(L-K)}{2}+1,$ is negative. Of course, $h^0$ is at least $1$, but it is not very likely in this case that it is greater than $1$. So most probably $|L|$ is not free. However, in the following example the estimate is $3$. I do not know if this example is the simplest one. Certainly, many other examples exist.

{\bf Example.} The following picture describes the graph of curves on a surface $X_1$, a hypothetical counterexample to the Jacobian Conjecture. Here the curve $A$ is of type $4$, the curve $B$ is of type $3$, the curves $D,$ $F_1,$ $F_2,$ $E,$ $F_1,$ $F_2,$ $F_3,$ and $F_4$ are of type $1$. All other curves are of type $2$. Curve $C$ is $\pi^{-1}(\infty)$.  Each vertex of $\Gamma$ is labeled by a $4$-tuple of numbers. The first three are the self-intersection, the $\bar{K}$-label, and the coefficient of $L$. The last one is the coefficient of $\Delta ,$ defined as follows.

Besides $|L|,$ the surface $X_1$ must posess another divisor class with interesting properties: the pullback by $\phi$ of the generic line on $Y$ that passes through the point $\phi (\pi^{-1}(\infty))$. This point is the image of all curves of type $2$, and this class of curves is $L-\Delta,$ where $\Delta =\sum d_iE_i,$ where $E_i$ are of type $2$. This class $\Delta$ has the following properties:

1) All coefficients $d_i$ are positive integers for curves of type $2$ (and zero for all other curves).

2) $\Delta$ intersects non-positively with all curves of type $2$.

3) $\Delta$ intersects by $1$ with all curves of type $1$ (which means that the coefficients $d_i$ adjacent to them are equal to $1$).

4) The intersection of $\Delta$ with any curve of type $3$ is less than or equal to that of $L.$

The above conditions are equivalent to a fairly complicated system of linear inequalities on the coefficients of $\Delta$ . However, for relatively small graphs like the one above, one can solve it using case-by-case analysis, paying attention to the "slope", which is the ratio $(d_i/a_i)$. This slope, as a function of the vertices, can not have a local minimum on any curve of type $2$.

\vskip 6cm

\includegraphics[viewport=0in 0in 5in 5in,scale=.2]{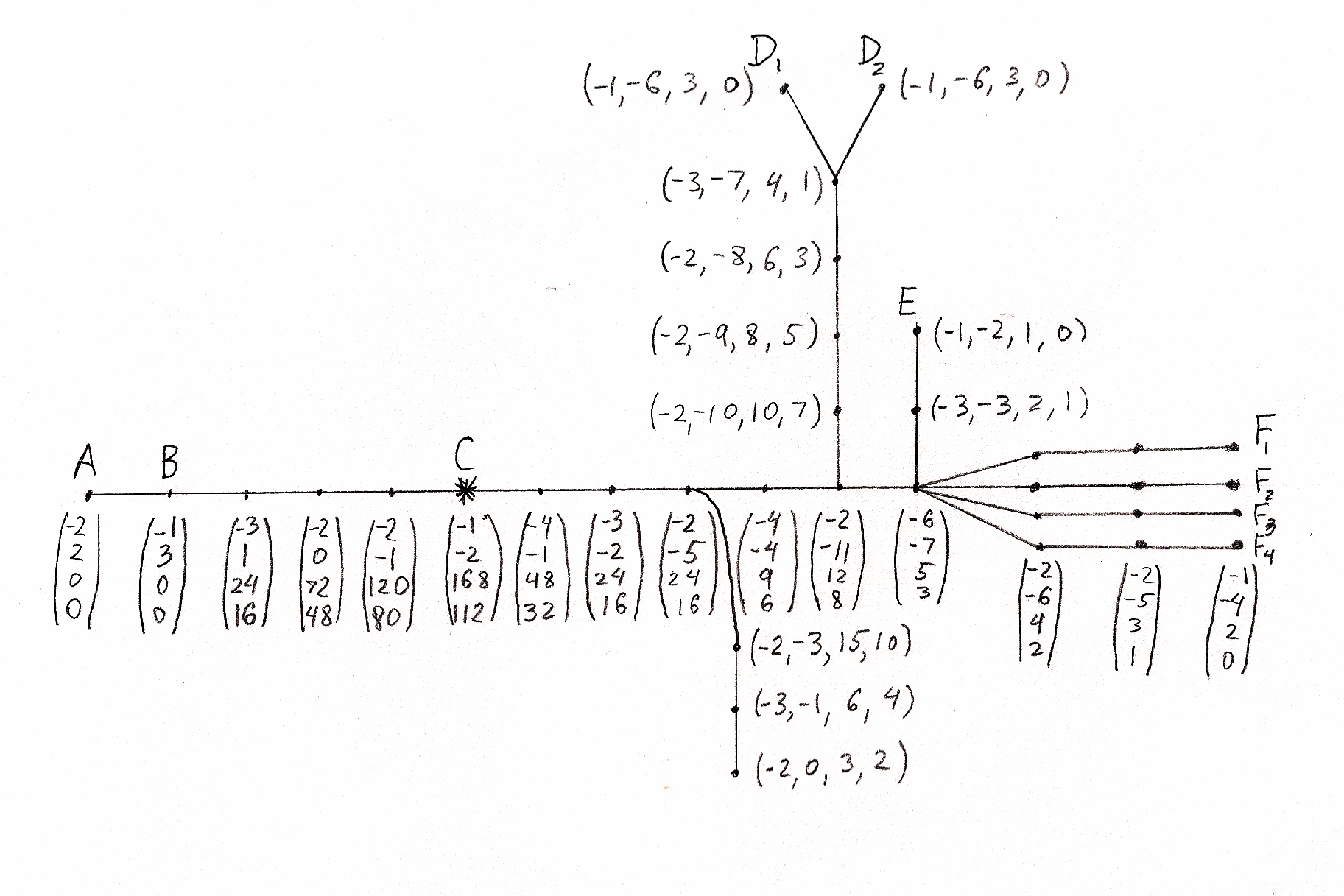}

The graph $\Gamma$ does not completely determine $X$ because there is a choice involved every time we blow up a  point on one of the curves. One can explicitly describe the moduli space of surfaces with the given graph as a Zariski open subset of some finite-dimensional affine space of parameters. If for some values of these parameters we happen to get $h^0(L)>h^0(L'),$ then we can get a map from $X$ to $P^2,$ by chosing a linear subspace of $H^0$ spanned by three sections: a section which has as its divisor the sum of the exceptional curves, $i(f)$, where $f$ is a generic section of $ L'$, and  a generic section of $H^0(L)$. This corresponds to polynomial map, given by two polynomials, of degrees $168$ and $56.$ This map would have the prescribed ramification at all curves of type $1$. It is possible, though by no means automatic that one can arrange for it to have the ramification of degree $3$ on the curve B. If this happens, it would be a counterexample to the Jacobian Conjecture.

The above example is the simplest one known to the author in which $(L^2-L\cdot K)/2 +1>0.$ One can certainly use  computers to systematically search for graphs like the example above, that can possibly provide a counterexample to the Jacobian Conjecture. Then one can try to solve the systems of linear  equations (with parameters) that determine $H^0(L)$ and $H^0(L-\Delta)$, and, finally, search for a counterexample to the Jacobian Conejcture.

If one believes that the Jacobian Conjecture is true in dimension two, one can try to push this research forward in several different ways. In most interesting examples the majority of curves are of type $2$. One can differentiate between them by blowing up points on the target $P^2$ (which will require further blowups on $X$). Alternatively, one can investigate what happens to the graph when the map $\phi$ is composed with an automorphism of the target $P^2$. This may lead to a new proof and possibly a strengthening of the results of Domrina and Orevkov (\cite{Domrina}, \cite{DomrinaOrevkov}) that the Jacobian Conjecture is true for maps of generic degree at most four.

If one blows up a sequence of points on $P^2$ so that the line at infinity is mapped to a curve, rather than to a point, one can contract all other exceptional curves on the target surface and get a rational surface with Picard number one with very special properties. It has a Zariski open subset isomorphic to $A^2.$ The complement is an irreducible rational curve $L$, with $L^2=\frac{1}{n}$ where $n$ is some natural number. The surface $X$ is singular, with rational singularities. It would be interesting to classify all such surfaces. The Jacobian Conjecture is directly related to some questions regarding the effective divisor classes on these surfaces.

\end{document}